\documentclass[11pt]{article}
\usepackage{latexsym,amssymb,amsmath,euscript,enumitem}
\usepackage{booktabs,color,pgf,tikz,pgfbaseimage}
\usepackage{mathrsfs}
\usetikzlibrary{arrows}
\usepackage{hyperref}

\title{On dendrites, generated by polyhedral systems and their ramification points.}

\author{Mary Samuel \and Andrey Tetenov \footnote{Supported by Russian Foundation of Basic Research project 16-01-00414}\and Dmitry Vaulin }

\begin{document}

\newcommand{\rr}{\mathbb{R}}
\newcommand {\nn} {\mathbb{N}}
\newcommand {\zz} {\mathbb{Z}}
\newcommand {\bbc} {\mathbb{C}}
\newcommand {\rd} {\mathbb{R}^d}
\newcommand {\rpo}{\mathbb{R}_+^1}

 \newcommand {\al} {\alpha}
\newcommand {\be} {\beta}
\newcommand {\da} {\delta}
\newcommand {\Da} {\Delta}
\newcommand {\ga} {\gamma}
\newcommand {\Ga} {\Gamma}
\newcommand {\la} {\lambda}
\newcommand {\La} {\Lambda}
\newcommand{\om}{\omega}
\newcommand{\Om}{\Omega}
\newcommand {\sa} {\sigma}
\newcommand {\Sa} {\Sigma}
\newcommand {\te} {\theta}
\newcommand {\fy} {\varphi}
\newcommand {\Fy} {\varPhi}

\newcommand {\rh} {\varrho}

\newcommand{\ep}{\varepsilon}

\newcommand{\VEC}{\overrightarrow}
\newcommand {\ra} {\rightarrow}
\newcommand{\IN}{{\subset}}
\newcommand{\ov}{{\overline}}
\newcommand{\NI}{{\supset}}
\newcommand \dd  {\partial}
\newcommand {\mmm}{{\setminus}}
\newcommand{\probel}{\vspace{.5cm}}
\newcommand{\8}{{\infty}}
\newcommand{\io}{{I^\infty}}
\newcommand{\ia}{{I^*}}
\newcommand{\0}{{\varnothing}}
\newcommand{\vse}{$\blacksquare$}

\newcommand {\bfep} {{{\bar \varepsilon}}}
\newcommand {\Dl} {\Delta}
\newcommand{\vA}{{\vec {A}}}
\newcommand{\vB}{{\vec {B}}}
\newcommand{\vF}{{\vec {F}}}
\newcommand{\vf}{{\vec {f}}}
\newcommand{\vh}{{\vec {h}}}
\newcommand{\vJ}{{\vec {J}}}
\newcommand{\vK}{{\vec {K}}}
\newcommand{\vP}{{\vec {P}}}
\newcommand{\vX}{{\vec {X}}}
\newcommand{\vY}{{\vec {Y}}}
\newcommand{\vZ}{{\vec {Z}}}
\newcommand{\vx}{{\vec {x}}}
\newcommand{\va}{{\vec {a}}}
\newcommand{\vga}{{\vec {\gamma}}}

\newcommand{\hf}{{\hat {f}}}
\newcommand{\hg}{{\hat {g}}}

\newcommand{\bj}{{\bf {j}}}

\newcommand{\bi}{{\bf {i}}}
\newcommand{\bk}{{\bf {k}}}
\newcommand{\bu}{{\bf {u}}}

\newcommand{\bX}{{\bf {X}}}
\newcommand{\Ep}{{\mathfrak p}}
\newcommand{\Eq}{{\mathfrak q}}
\newcommand{\eJ}{{\EuScript J}}
\newcommand{\wP}{{\widetilde P}}
\newcommand{\eU}{{\EuScript U}}
\newcommand{\eS}{{\EuScript S}}
\newcommand{\eH}{{\EuScript H}}
\newcommand{\eC}{{\EuScript C}}
\newcommand{\eP}{{\EuScript P}}
\newcommand{\eT}{{\EuScript T}}
\newcommand{\eG}{{\EuScript G}}
\newcommand{\eK}{{\EuScript K}}
\newcommand{\eF}{{\EuScript F}}
\newcommand{\eZ}{{\EuScript Z}}
\newcommand{\eL}{{\EuScript L}}
\newcommand{\eD}{{\EuScript D}}
\newcommand{\E}{{\EuScript E}}
\def \diam {\mathop{\rm diam}\nolimits}
\def \sup {\mathop{\rm sup}\nolimits}
\def \fix {\mathop{\rm fix}\nolimits}
\def \Lip {\mathop{\rm Lip}\nolimits}
\def \min {\mathop{\rm min}\nolimits}

\newcommand{\red}{\textcolor{red}}

\newtheorem{thm}{\bf Theorem}
 \newtheorem{sled}[thm]{\bf Corollary}
 \newtheorem{lem}[thm]{\bf Lemma}
 \newtheorem{prop}[thm]{\bf Proposition}
 \newtheorem{dfn}[thm]{\bf Definition}
 \newcommand{\rmk}{{\bf {Remark.}}}
\newcommand{\dok}{{\bf{Proof:  }}}

\maketitle

\bigskip

\begin{abstract}The paper considers   systems of contraction similarities in $\mathbb R^d$ sending a given polyhedron $P$ to polyhedra $P_i\subset P$, whose non-empty intersections  are singletons and contain the common vertices of those polyhedra, while the intersection hypergraph of the system is acyclic. It is proved that the attractor $K$ of such system is a dendrite  in $\mathbb R^d$. The ramification points of such dendrite fave finite order whose upper bound depends only on the polyhedron $P$, and the set of the cut points of  the dendrite $K$ is equal to the dimension of the whole $K$ iff $K$ is a Jordan arc.
\end{abstract}

\section{Introduction}

Though the study of topological properties of dendrites from the viewpoint of general topology proceed for more than  three quarters of a century \cite{Char,Kur, Nad}, the attempts to study the geometrical properties of self-similar dendrites are rather fragmentary.

In 1985 M.~Hata \cite{Hata}  studied the connectedness properties of self-similar sets and proved that if a dendrite is an attractor of a system of weak contractions in a complete metric space, then the set of its enpoints is infinite. In 1990 Ch.~Bandt showed in his unpublished paper \cite{BS} that the Jordan arcs connecting pairs of points of a post-critically finite self-similar dendrite are self-similar, and the set of possible values for dimensions of such arcs is finite. Jun Kigami in his work \cite{Kig95} applied the methods of harmonic calculus on fractals to dendrites; on a way to this he developed effective approaches to the study of  structure of self-similar dendrites. D.Croydon in his thesis \cite{C} obtained heat kernel estimates for continuum random tree and for certain family of p.c.f. random dendrites on the plane. D.Dumitru and A.Mihail \cite{DM} made an attempt to get a sufficient condition for a self-similar set to be a dendrite in terms of sequences of  intersection graphs for the refinements of the system $\eS$.

There are many papers \cite{Bar,BK,Zell} discussing examples of self-similar dendrites, but systematic approach to the study of self-similar requires to find the answer to the following questions:
What kind of topological restrictions characterise the class of dendrites generated by systems of similarities in $\rr^d$? What are the explicit construction algorithms for self-similar dendrites? What are the metric and analytic properties of morphisms of self-similar structures on dendrites?

To approach these questions, we start from simplest and  most obvious settings, which were used by many authors \cite{BS,Str}.   We consider systems $\eS$ of contraction similarities in $\rd$
defined by some polyhedron $P\IN\rd$, which we call contractible   $P$-polyhedral systems.

We prove that the attractor of such system $\eS$ is a dendrite $K$ in $\rd$ (Theorem \ref{main}), and there is a dense subset of $K$ such that punctured neigh\-bour\-hoods of its  points split to a finite disjoint union of subsets of solid angles  $\Om_l$, equal to the solid angles of $P$ (Theorem \ref{refsys}); we show that the orders of points $x\in K$ have an upper bound, depending only on $P$ (Theorem \ref{order}); and that
Hausdorff dimension of the set $CP(K)$ of the cut points of $K$ is strictly smaller than the dimension of the set $EP(K)$ of its end points unless $K$ is a Jordan arc (Theorem \ref{equal}).

\subsection{Preliminaries}

{\bf Dendrites.}  A {\em dendrite} is a locally connected continuum containing no simple closed curve.

In the case of dendrites the order $Ord(p,X)$ of the point  $p$  with respect to  $X$ is equal to the number  of components of the set $X \setminus \{p\}$. the points of order
 $1$ are called {\em end points} in $X$, and cut points are called usual points if  $Ord(p,X)=2$ and {\em ramification points}, if $Ord(p,X)\ge
3$.

We will  use the following statements selected from  \cite[Theorem 1.1]{Char}:
\begin{thm}\label{Char} For a continuum $X$ the following conditions are equivalent:
\begin{enumerate}[label=(\alph*),noitemsep]
\item $X$ is dendrite;
\item every two distinct points of $X$ are separated by a third point;
\item each point of $X$ is either a cut point or an end point of $X$;
\item each nondegenerate subcontinuum of $X$ contains uncountably many cut points of $X$.
\item for each point $p \in X$ the number of components of the set $X \setminus \{p\} = ord (p, X)$ whenever either of these is finite;
\item the intersection of every two connected subsets of X is connected;
\item X is locally connected and uniquely arcwise connected.
\end{enumerate}
\end{thm}

\bigskip
{\bf Self-similar sets.}
Let $(X, d)$ be a complete metric space. 
A mapping\\ $F: X \to X$ is a contraction if $\Lip F < 1$. 
The mapping $S: X \ra X$ is called a similarity if \begin{equation} d(S(x), S(y)) = r d(x, y) \end{equation} for all $x, y\in X$ and some fixed r. 

\begin{dfn} 
Let $\eS=\{S_1, S_2, \ldots, S_m\}$ be a system of (injective) contraction maps on the complete metric space $(X, d)$.
 A nonempty compact set $K\IN X$ is said to be invariant with respect to $\eS$, if $K = \bigcup \limits_{i = 1}^m S_i (K)$. \end{dfn}
 
 We also call the   subset $K \IN X$ self-similar with respect to $\eS$.\\
Throughout the whole paper, the maps $S_i\in \eS$ are supposed to be  similarities and the set $X$ to be $R^d$.\probel

{\bf Notation.} $I=\{1,2,...,m\}$ is the set of indices, $\ia=\bigcup\limits_{n=1}^\8 I^n$ -
is the set of all finite $I$-tuples, or multiindices $\bj=j_1j_2...j_n$,  where $\bi\bj$  is the concatenation of the corresponding multiindices;\\ 
we say $\bi\sqsubset\bj$, if $\bi=i_1\ldots i_n$ is the initial segment in $\bj=j_1\ldots j_{n+k}$ or $\bj=\bi\bk$ for some $\bk\in\ia$;
if $\bi\not\sqsubset\bj$ and $\bj\not\sqsubset\bi$, $\bi$ and $\bj$ are {\em incomparable};\\
we write
$S_\bj=S_{j_1j_2...j_n}=S_{j_1}S_{j_2}...S_{j_n}$ 
and for the set $A
\subset X$ we denote $S_\bj(A)$ by $A_\bj$; 
we also denote by $G_\eS=\{S_\bj, \bj\in\ia\}$ the semigroup, generated by $\eS$;\\
$I^{\8}=\{{\bf \al}=\al_1\al_2\ldots,\ \ \al_i\in I\}$ --
index space; and $\pi:I^{\8}\rightarrow K$ is the {\em index map
 }, which sends $\bf\al$ to  the point $\bigcap\limits_{n=1}^\8 K_{\al_1\ldots\al_n}$.\\

\begin{dfn}The system ${\eS}$ satisfies the {\em open set condition} (OSC) if there exists a non-empty open set $O \IN X$ such that $S_i (O), \{1 \le i\le m\}$ are pairwise disjoint and all contained in $O$.\end{dfn}

We say the self-similar set $K$ defined  by  the  system $\eS$  satisfies the one-point intersection property if for  any $i\neq j$,  $S_i(K) \bigcap S_j(K)$ is not more than one point. 
\\

We use  the following  convenient criterion of  connectedness of  the  attractor  of a system $\eS$    \cite{Hata,Kig}.

\begin{thm}\label{Kig} Let $K$ be the attractor of a system of contractions $\eS$  in a complete metric space $(X,d)$. Then the following statements are equivalent:\\
1) for any $i, j \in I$ there are $\{i_0, i_1, \ldots, i_n\} \IN I$ such that $i_0 = i, i_n = j$ and \\$S_{i_k}(K)\bigcap S_{i_{k + 1}}(K) \ne \0$ for any $k = 0, 1, \ldots, {n-1}$.\\
2) $K$ is arcwise connected.\\
3) $K$ is connected.\end{thm}

\begin{prop}
If a self-similar set $K$ is connected, it is locally connected.
\end{prop}

{\bf Zippers and multizippers.} The simplest  way to construct  a self-similar  curve is  to  take a polygonal line and then replace each of its segments by a smaller copy of the same polygonal line; this construction is called   zipper and was studied  by Aseev, Tetenov and  Kravchenko \cite{ATK}.

\begin{dfn} Let $X$  be a complete metric space. A system $\eS = \{S_1, \ldots, S_m\}$
of contraction mappings of $X$ to itself is called a {\em zipper} with vertices\\ $\{z_0, \ldots, z_m\}$
and signature $\ep = (\ep_1, \ldots, \ep_m), \ep_i \in\{0,1\}$, if for    $ i = 1\ldots m$, $S_i (z_0) = z_{i-1+\ep_i}$ and $S_i (z_m) = z_{i-{\ep_i}}$.\end{dfn}

More  general approach  for  building self-similar  curves and  continua is provided  by  a graph-directed version of zipper construction \cite{Tet06}:

\begin{dfn} Let $\{X_u, u \in V\}$ be a system of spaces, all isomorphic to ${\rr}^d$.
 For each $X_u$ let a finite array of points be given $\{x_0^{(u)}, \ldots, x_{m_u}^{(u)}\}$.
 Suppose for each $u\in V$ and $0\le k \le m_u$ we have some $v (u, k) \in V$  and  $\ep(u,k)\in\{0,1\} $   and a map ${\eS}_k^{(u)}: X_v \ra X_u$ such that\\
 $S^{(u)}_k (x_0^{(v)}) = x_{k-1}^{(u)}$  or $x_k^{(u)}$  and  
$S_k^{(u)} (x_{m_v}^{(v)}) = x_{k}^{(u)}$  or  $  x_{k-1}^{(u)}$,
depending on the signature $\ep (u, r)$.\\ The graph directed iterated function system (IFS) defined by the maps $S_k^{(u)}$ is called a {\em multizipper} ${\eZ}.$\end{dfn}

The attractor of multizipper ${\eZ}$ is a system of  connected and  arcwise connected  compact  sets $K_u\IN X_u$ satisfying  the system of equations
$$ K_u=\bigcup\limits_{k=1}^{m_u}S_k^{(u)}(K_{v(u,k)}),\qquad u\in V$$
 We call the  sets $K_u$  the components of  the attractor of ${\eZ}$.

The components $K_u$  of the  attractor of $\eZ$ are  Jordan  arcs if  the following  conditions  are  satisfied:
 
\begin{thm} \label{Jor}
 Let ${\eZ}_0 = \{S_k^{(u)}\}$ be a multizipper with node points $x_k ^{(u)}$ and a signature $\ep = \{(v (u, k), \ep (u, k)), u \in V, k = 1, \ldots, m_u\}$. If for any $u \in V$
and any  $ i, j \in \{1, 2, \ldots, m_u\}$, the set $K_{(u, i)} \cap K_{(u, j)} = \0$ if $|i- j| > 1$ and is a singleton if $|i -j| = 1$, then any linear parametrization $\{f_u: I_u \ra K_u\}$ is a homeomorphism and each $K_u$ is a Jordan arc with endpoints $x^{(u)}_0,\   x_m ^{(u)}.$\end{thm}

\section{Contractible polyhedral systems.}
Let $P\IN \rd$ be a finite polyhedron  homeomorphic to a  $d$-dimensional ball  and let $V_P=\{A_1,...,A_{n_P}\}$ be the set of its vertices, and $\Om(P,A_i)$ be the solid angles at the vertices of $P$.\smallskip\\
Consider a system of similarities  $\eS = \{S_1, \ldots, S_m\}$, which define polyhedra $P_i=S_i(P)$ and satisfy the following conditions:\\
{\bf(D1)}\ \  For any $i\in I$,  $P_i\IN P$;
\\
{\bf(D2)}\ \  For any $i\neq j,\ \   i, j \in I,$ the intersection $P_i\cap P_j$ is either empty, or is a common vertex of   $P_i$ and$P_j$;\\
{\bf(D3)}\ \  $V_P\IN \bigcup\limits_{i\in I}S_i(V_P)$;\\
{\bf(D4)}\ \   { The set    ${\wP} = \bigcup \limits_{i = 1}^m P_i$ is contractible.}

\begin{dfn}\label{pts}
A system $\eS $, satisfying {\bf D1-D4},
is called  $P$-polyhedral system of similarities.
\end{dfn}

The similarities $S_i\in\eS$ are contractions, therefore the system $\eS$ has the attractor $K$; the system  $\eS$ поgenerates the semigroup $G_\eS=\{S_\bj,\bj\in \ia\}$ and therefore defines the set of polyhedra  $G_\eS(P)=\{P_\bj, \bj\in\ia\}$. The properties of this system of refining polyhedra define the geometric properties of the invariant set $K$. First we focus on those properties, which follow from {\bf D1--- D3} only, which corresponds to a class of point connected self-similar sets, as they were defined by R.Strichartz \cite{Str}.
The reative position of solid angles of polyhedra $P_\bj$ will be our special interest:

\begin{thm}\label{refsys} Let  $\eS$ be a $P$-polyhedral  system  of similarities.\\
(a) The system $\eS$ satisfies (OSC).\\
(b) $P_\bj\IN P_\bi$ iff  $\bj\sqsupset\bi$.\\
(c) If $\bi\sqsubset\bj$, then  $ S_\bi(V_P)\cap P_\bj\IN S_\bj(V_P)$.\\
(d) For any incomparable $\bi,\bj\in \ia$, $\#(P_\bi\cap P_\bj)\le 1$ and $P_\bi\cap P_\bj=S_\bi(V_P)\cap S_\bj(V_P)$.\\
(e) The set $G_\eS(V_P)$ of vertices of polyhedra $P_\bj$ is contained in $K$.\\
(f)If $x\in K\mmm G_\eS(V_P)$, then $\#\pi^{-1}(x)=1$.\\
(g) For any $x\in G_\eS(V_P)$ there is   $\ep>0$ 
and a finite system $\{\Om_1,...,\Om_n\},$ where $n=\#\pi^{-1}(x)$,
of disjoint solid angles with vertex  $x$,
such that if $x\in P_\bj$ and $\diam P_\bj<\ep$, then for some
  $k\le n$, $\Om(P_\bj,x)=\Om_k$. Conversely, for any $\Om_k$ there is such $\bj\in\ia$, that  $\Om(P_\bj,x)=\Om_k$.
\end{thm} 
\dok (a)  It follows from {\bf D1, D2} that the interior of $P$   is the desired  open set for (OSC);  (b) follows from (OSC);
 
 (c) Notice that {\bf D2, D3} imply the condition   {\bf (D3a)}: for any $i\in I$, $ P_i\cap V_P\IN S_i(V_P)$:\\ Indeed, if $x\in P\mmm V_p$ and $S_i(x)=A\in V_P$,  then, since there is  $j\in I$, such that $A\in S_j(V_P)$, $P_i\cap P_j\notin S_i(V_P)$, whis=ch contradicts {\bf D3}. 
 
 Using induction, we derive from {\bf D3a }, that for any $\bk \in \ia$, $ P_\bk\cap V_P\IN S_\bk(V_P)$;
 
 Let now $\bj=\bi\bk$ and $A\in S_\bi(V_P)\cap S_\bi(P_\bk)$. It means that $S_\bi^{-1}(A)\in V_P\cap P_\bk$, and therefore $S_\bi^{-1}(A)\in S_\bk(V_P)$, or
$A\in S_\bj(V_P)$.
 
    (d) Represent a pair of incomparable multiindices as $\bk\bi,\bk\bj$, where $i_1\neq j_1$. Since
 $P_{\bk\bi}\cap P_{\bk\bj}\neq\0$,   
 $P_\bi\cap P_\bj\neq\0$.  
    But $P_\bi\cap P_\bj\IN P_{i_1}\cap P_{j_1}$. 
    The last set is nonempty and therefore it consists solely of a common vertex of $P_{i_1}$ and $P_{j_1}$; 
     by (c), this point is is also a common vertex of $P_\bi$ and$P_\bj$; 
    therefore $P_{\bk\bi}\cap P_{\bk\bj}=S_{\bk\bi}(V_P)\cap S_{\bk\bj}(V_P)$.
    
    (e) For any vertex $A\in V_P$ there are $A_1\in V_P$  and
     $\al_1\in I$ such that $S_{\al_1}(A_1)=A$. By induction we get that for any $n$  there are such $A_n\in V_P$ and  $\al_1\ldots\al_n\in I^n$, 
     that $S_{\al_1\ldots\al_n}(A_n)=A$.
      In this case, $\bigcap\limits_{n=1}^\8 S_{\al_1\ldots\al_n}(P)=\{A\}$ and $A\in K$.
       Thus, $V_P\IN K$, and therefore $G_\eS(V_P)\IN K$.
       
 (f)  If $\pi^{-1}(x)$ contains non-equal ${\bf \al,\be}\in I^\8$, then for some  $n$, $\al_1\ldots\al_n$ and $\be_1\ldots\be_n$ are incomparable; therefore $x\in P_{\al_1\ldots\al_n}\cap P_{\be_1\ldots\be_n}$, so $x\in G_\eS(V_P)$.
 
 (g) First let ${\bf\al}\in \io$ and $\pi( {\bf\al})=A\in V_P$. 
 As in (e), for any $n$, $S_{\al_1\ldots\al_n}(A_n)=A$ and 
  $S_{\al_1\ldots\al_n}(\Om(P,A_n))\IN \Om(P,A)$. Moreover, the solid angles   $S_{\al_1\ldots\al_n}(\Om(P,A_n))$  form a nested sequence.  Since the set $\{\Om(P,B), B\in V_P\}$ is finite, there is a solid angle $\Om_\al$ and a number $N\in\nn$,
  such that if $n>N$, then $S_{\al_1\ldots\al_n}(\Om(P,A_n))= \Om_\al$. 
  At the  same time,  $S_{\al_1\ldots\al_n}(P)\IN \Om_\al$. If for some
  ${\bf\be}\in \io, {\bf\be}\neq {\bf\al}$, $\pi( {\bf\be})=A$, then,
   according to (d), $\Om_\al\cap \Om_\be=\{A\}$. Thus, the set $\pi^{-1}(A)$ can be mapped bijectively to the family of disjoint solid angles  $\Om_k$ with common vertex $A$.

  The measure
 $\te(\Om_k)$
 is greater or equal to
  $\te_{min}=$ $\min\{\te(\Om(P,A)), A\in V_P\}$, therefore the set of different $ {\bf\al}\in \io$ such that $\pi({\bf\al})=A$, does not exceed
  $\te(\Om(P,A))/\te_{min}$, if $A\in V_P$, and  $\te_F/\te_{min}$, if $A\in\dot P$, where $\te_F$ -- is the measure of complete solid angle in $\rd$. 
   \vse
   
 Now we discuss some properties of $\eS$ which follow from the condition {\bf D4}. 

 Applying Hutchinson operator  $T (A)=\bigcup\limits_{i\in I} S_i(A)$ of the system $\eS$ to the polyhedron $P$, we get the set ${\wP} = \bigcup \limits_{i\in I} P_i$. We define ${\wP}^{(1)} = T(P), {\wP}^{(n + 1)} = T(\wP^{(n)})$. Thus we get a nested family of sets
${\wP}^{(1)}\NI {\wP}^{(2)}\NI \ldots \NI  {\wP}^{(n)}\NI\ldots,$
whose intersection is  $K$.

\bigskip
The composition of two contractible
$P$-polyhedral systems is also of the same type:

\begin{lem}\label{trans} Let $\eS$ and ${\eS'}$ be contractible $P$-polyhedral systems of similarities.
Then  ${\eS}'' = \{S_i \circ S_j', S_i \in {\eS}, S_j \in
{\eS}'\}$ --- is also contractible $P$-polyhedral system of similarities.\end{lem}

\noindent{\dok} {\bf(D1)}\ \ is obvious since $S_i \circ S_j'(P) \IN S_i (P) \IN P$.

{\bf(D2)}\ \ Let $Q_1 =S_{i_1}\circ S_{j_1}' (P)$ and $Q_2 = S_{i_2}\circ S_{j_2}' (P)$ be two polyhedra for the system ${\eS}''$; consider their intersection:\\
if $i_1 \ne i_2$, then $Q_1 \bigcap Q_2 \IN P_{i_1} \bigcap P_{i_2}$,
  where the right part is either empty, or for some $A_1, A_2\in V_P$,  $P_{i_1}\cap P_{i_2}=\{S_{i_1}(A_1)\}= \{S_{i_2}(A_2)\}$. Since $A_1\in S'_{j_1}(V_P)$ и $A_2\in S_{j_2}'(V_P)$, $Q_1\cap Q_2 =S_{i_1}\circ S_{j_1}' (V_P) \cap S_{i_2}\circ S_{j_2}' (V_P)$;\\     
    if $i_1 = i_2$, then  $Q_1 \bigcap Q_2 = S_{i_1} (P_{j_1}' \bigcap P_{j_2}')$
 where the right part is either empty or a  singleton contained in $S_{j_1}'(V_P)\cap S_{j_2}'(V_P)$.

{\bf(D3)}\ \   holds, because for any vertex $A\in V_P$, there are $A_1\in V_P$ and $S_{i_1}\in \eS$ such that $S_{i_1}(A_{1})=A $; ,also there are  $S_{i_2}'\in \eS'$ and $A_{2}\in V_P$ such that$S_{i_2}'(A_2)=A_1$; therefore $S_{i_1}S_{i_2}'(A_2)=A$. If $x\in P$ and $S_{i_1}S_{i_2}'(x)=A$, then $S_{i_2}'(x)\in V_P$, therefore $x\in V_P$.

{\bf(D4)}\ \ The sets $\wP=\bigcup \limits_{i = 1}^m P_i$ and
$\wP'=\bigcup \limits_{i = 1}^{m'} P_i'$ are strong deformation retracts of  $P$, containing the set $V_P$. Let ${\fy}' (X, t): P \times [0, 1] \ra P$ be the deformation retraction
from  $P$ to $\mathop{\bigcup}\limits_{i =
1}^{m'} P_i'$. The map $\fy'$ satisfies the following conditions:
 ${\fy}' (x, 0) = Id$, ${\fy}' (x, 1) (P) = {\wP}'$ and for any $t\in[0,1]$,
${{\fy}' (x, t)}|_{{\wP}'} = {Id}_{{\wP}'}$.

Define the map ${\fy}_i' : {P_i} \times [0, 1] \ra P_i$  
by a formula
$${\fy}_i' (x, t) = S_i \circ {\fy}' ({S_i}^{-1} (x), t).$$ Each map
${\fy}_i'$ is a deformation retraction from $P_i$ to $S_i
({\wP}')$.
\\
Observe that the map ${\fy}_i'$ keeps all the vertices $S_i (A_k)$
of the polyhedron $P_i$ fixed. Therefore we can define a deformation
retraction ${\widetilde\fy} (x, t):{\wP} \times [0, 1] \to
 \bigcup\limits_{i = 1}^m{S_i ({\wP}')} ={\wP}$  by a formula
$${\widetilde\fy} (x, t) =  \fy_i' (x, t),
\mbox{\quad \rm  if  }x \in P_i$$
 The map $\widetilde\fy$ is well-defined and continuous
because if $P_i\bigcap P_j = \{S_i (A_k)\}=\{S_j(A_l)\}$ for some
$k$ and   $l$, then
${\fy}_i' (S_i (A_k), t) \equiv {\fy}_j' (S_j (A_l), t) \equiv S_i (A_k)$.
\\
 Moreover,${\widetilde\fy} (x, 0) =x$ on $\wP$, and
${\widetilde\fy} ({\wP}, 1) \equiv \bigcup\limits_{i = 1}^m S_i ({\wP}')$ and ${\widetilde\fy} (x, t)|_ {{\wP}''} \equiv Id$. \\ So ${\widetilde\fy} (x, t)$ is a deformation retraction from $\wP$ to ${\wP}''$.\\
Therefore, the set ${\wP}'' = \bigcup {S_i \circ S_j'(P)}$ is contractible.{\vse}
\begin{sled}\label{sntoo} If $\eS$ is a contractible $P$-polyhedral system, the same is true for $\eS^{(n)}=\{S_\bj,\bj\in I^n\}$. 
\vse
\end{sled}

The contractibility of the set $\wP$ and the condition {\bf D2} imply, that any simple closed curve in $\wP$
lies in one of the polyhedra $P_i$; this can be derived from the following  Lemma:

\begin{lem}\label{balls} Let $B_i, i=1,...,n$ --  be a finite family of topological balls, such that for any $i,j$ the intersection $B_i\cap B_j$ contains no more than one point and the set $X=\bigcup\limits_{i=1}^n B_i$ is simply-connected. Then any simple closed curve in $X$ lies in some $B_i$.\end{lem}

\dok Choose in each of $B_i$ a point $O_i\in\dot B_i$ and for each $\{p_{ij}\}=B_i\cap B_j$ take a Jordan arc
  $\ga_{ij}$ with endpoints $O_i$ and $p_{ij}$ so that $\ga_{ij}\cap\ga_{ij'}=\{O_i\}$ if $j'\neq j$. Let $\Ga$ be a topological graph with vertices $O_i, i=1,...,n$  and $p_{ij}$ whose edges are $\ga_{ij}$. 
 Since for any  $i$ the union $\bigcup\limits_{j} \ga_{ij}$ is a strong deformation retract of the ball $B_i$, $\Ga$ is a strong deformation retract of the set  $X$ and therefore the graph $\Ga$ is a tree.
 
 Let $l$ be some Jordan arc in $X$. Suppose $l$ is in general position in the sence that $p_{ij}\in l$ iff $l\cap \dot{B_i}\neq\0$ and $l\cap \dot{B_j}\neq\0$.
 Each point $p_{ij}$ splits  $X$ to no less than 2 components. Therefore if $l\ni p_{ij}$, the arc $l$ is not closed.
 Thus, any simple closed curve in $X$ lies completely in one of the balls $B_i$.    \vse

\begin{thm}\label{main} The attractor $K$ of contractible $P$-polyhedral system of si\-mi\-la\-ri\-ties $\eS$ is a dendrite.\end{thm}

{\dok}   By Corollary \ref{sntoo}, the sets
${\wP}^{(n)}$ are contractible, compact and and satisfy
${\wP}^{(1)} \NI {\wP}^{(2)}\NI  {\wP}^{(3)} \ldots$. The diameter of components of the interior of any of ${\wP}^{(n)}$ does not exceed  ${\diam}{P}\cdot {q^n}$, where $q = \max\Lip(S_i)$.
Thus the set $K = \bigcap {\wP}^{(n)}$ is connected and has empty interior. Since $K$ is connected, it is locally connected and arcwise connected \cite[Theorem 1.6.2, Proposition
1.6.4]{Kig}. 

Let $l$ be some   Jordan arc in $K$. For any $n\in \nn$, $l\IN\wP^{(n)}$, so it follows from  Lemma \ref{balls} that if  $l$ has non-zero dianmeter, it is not closed. Therefore $K$ is a dendrite. {\vse}

The dendrite $K$ is contained in the polyhedron $P$; in general, its in\-ter\-sec\-tion with the boundary of  
 $P$ may be uncountable or  it  can contain even some whole edges of $P$. The same is also true for the intersection of the dendrite $K$ with each polyhedron
$S_\bj(P),\bj\in I^*$. Nevertheless it follows from {\bf D2}  that a subcontinuum $L\IN K$ can ''penetrate'' to a polyhedron  $S_\bj(P)$ only through its vertices, namely:

\begin{prop}\label{squeeze}
Let $\bj\in I^*$ be a multiindex. For any continuum $L\IN K$, whose intersection with both $ P_\bj$ and its exterior $\dot {CP_\bj}$ is nonempty, the set $\overline{L\mmm P_\bj}\cap P_\bj\IN S_\bj(V_P)$.
\end{prop}
\dok Observe that for any polyhedron $P_\bj, \bj\in
I^k$ the set ${\wP}^{(k)}\mmm S_\bj(V_P)$ is disconnected, and $P_\bj\mmm S_\bj(V_P)$ is its connected component, whose intersection with $K$ is equal to
$S_\bj(K\mmm S_\bj(V_P)$. Therefore   
$L\mmm S_\bj(V_P)$ is also disconnected.\vse

\subsection{The main tree and ramification points}
Since $K$ is a dendrite, by Theorem \ref{Char} for any vertices $A_i,A_j \in V_P$ there is unique Jordan arc ${\ga}_{ij}\IN K$ connecting $A_i,A_j$.
As it was proved by C.~Bandt \cite{BS}, these arcs are the components of the attractor of a graph-directed system of similarities. We show that this  system is a Jordan multizipper \cite{Tet06}:

\begin{thm}\label{arcs} The arcs ${\ga}_{ij}$ are the components of the attractor of some Jordan multizipper $\eZ$.\end{thm}
{\dok} We say, that polyhedra  $P_{i_1}, \ldots, P_{i_s}, i_k\in I$
form a chain, connecting  $x$ and $y$, if $x\in P_{i_1}, y\in P_{i_s}
$ and the intersection $ P_{i_k}\bigcap P_{i_l}$ is empty, if $|l - k| > 1$ and
is a common vertex of $ P_{i_k}$ and $ P_{i_l}$, if $| l - k| =
1$.
\\
For the vertices $A_i, A_j$, there is unique chain of polyhedra in
 the $P$-polyhedral system $\eS$, which consists of those   $P_k$, 
for which $\#P_k\cap\ga_{ij}\ge 2$; we  denote the polyhedra forming the chain and corresponding maps as   $P'_{ijk}=S'_{ijk}(P), k = 1, \ldots m_{ij}$, 
keeping in mind that all $S'_{ijk}\in\eS$.
\\
Let $u (i, j, k)$ и $v (i, j, k)$ be the indices of vertices $P$, for  which\\ $S'_{ijk}(A_u) = P'_{ij{(k-1)}} \bigcap P'_{ijk}= z_{ij{(k - 1)}}$\\
and $S'_{ijk}(A_v) =  P'_{ij{k}} \bigcap P'_{ij{(k+1)}}= z_{ijk}$, if $1 < k < m_{ij}$,\\ and if $k=1$ or $k =m_{ij}$, $u (i, j,1) = A_i= z_{ij0}$ and $v (i, j, m_{ij})= A_j = z_{ijm_{ij}}$. \\

Thus for any triple $(i,j,k)$,$1\le k\le m_{ij}$,   such $u,v\in\{1,...,n_P\}$ are specified, that
$S'_{ijk}(z_{uv0})=z_{ij(k-1)}$ and 
$S'_{ijk}(z_{uvm_{uv}})=z_{ijk}$.

Therefore the system  $\{S'_{ijk}\}$ is a multizipper $\eZ$ with node points $z_{ijk}$.\\

Since the relations: $${\ga}_{ij} = \bigcup\limits_{i = 1}^{m_{ij}} S'_{ijk}({\ga}_{u (i,j,k), v(i, j, k)}) =
\bigcup\limits_{i = 1}^{m_{ij}}{\ga}_{ijk} $$ 
are satisfied, the arcs  $\ga_{ij}$ form a complete set of the components of the attractor of the multizipper $\eZ$.

Since each arc   ${\ga}_{ijk}$ lies in  $P_{ijk}$, \\
$$ {\ga}_{ijk} \bigcap {\ga}_{ijl} = \0,$$ if $|k - l| > 1$ and $${\ga}_{ijk} \bigcap {\ga}_{ijl} =\{ z_{ijk}\},$$ and $l = k \pm 1$.  
Therefore $\eZ$ satisfies the conditions of the Theorem  \ref{Jor}  and is a Jordan  multizipper. {\vse}

The set  $\hat\ga=\bigcup\limits_{i\neq j}\ga_{ij}$ is a subcontinuum of the dendrite  $K$ and therefore is a dendrite. Since all the end points of $\hat\ga$  are contained in $V_P$,  $\hat\ga$ is a finite dendrite or topological tree  \cite[{\bf A.17}]{Char}. Let $n_E$ be the number of end points of  $\hat\ga$. As it was pointed out by Kigami \cite{Kig95}, $\hat\ga$   may be represented as  union of at most
 $(n_E-1)$ Jordan arcs having disjoint interiors.

\begin{dfn}\label{defmt}
The union $\hat\ga=\bigcup\limits_{i\neq j}\ga_{ij}$ is called 
{\em the main tree} of the dendrite $K$. The ramification points of
$\hat\ga$ are called {\em main ramification points} of the dendrite
$K$.
\end{dfn}


The following statement establishes the relations between the sets of vertices   $V_P$,    end points $EP(\hat\ga)$ and  cut points $CP(\hat\ga)$    of the main tree $\hat\ga$:

\begin{prop}\label{xintree} Let $x\in K$.\\
(a)  $\hat\ga~\IN\bigcup\limits_{ A_j\in V_P}\ga_{A_jx}$; besides, if $\hat\ga~=\bigcup\limits_{ A_j\in V_P}\ga_{A_jx}$, then $x\in\hat\ga$;\\
 (b) $ EP(\hat\ga)=V_P\mmm CP(\hat\ga)$;\\
 (c) $x\in CP(\hat\ga)$ iff there are vetrices $A_{i}$, $A_{j}$, belonging to different components of
$K\mmm\{x\}$;\smallskip \\
(d) for $x\in CP(K)$,   $Ord(x,K)=Ord(x,\hat\ga)$ iff for any component $C_l$ of the set $K\mmm\{x\}$, $C_l\cap V_P\neq\0$.
\end{prop}
\dok  For any $A_i,A_j\in V_P$, $\ga_{A_iA_j}\IN
\ga_{A_{i}x}\cup\ga_{A_{j}x}$, which gives (a).
To get (b), notice that if $x\in\hat\ga$ is not a vertex  then $x$ it is the inner point of some arc
 $\ga_{A_iA_j}$, therefore it is a cut point of  
$\hat\ga$ and therefore $x\notin EP(\hat\ga)$. \\
(c): Since $\ga_{A_ix}\cap \ga_{A_jx}=\{x\}$, we have
$\ga_{xA_{i}}\cup\ga_{xA_{j}}=\ga_{A_iA_j}$. So $x$ is a cut point of $\ga_{A_iA_j}$,
and therefore of $\hat\ga$. \\ (d): Necessity is obvious, so we prove sufficiency. By (c), $x\in CP(\hat\ga)$. The number of components of $K\mmm\{x\}$ is no greater than $n_P$, so  $Ord(x,K)$ is finite. Let
$C_l,l=1,...,k, k=Ord(x,K)$ be the components of $K\mmm\{x\}$. It also follows from (c) that $x\in\hat\ga$ and that the two vertices $A_{i }$ and $A_{j}$ lie in the same component $C_l$ iff 
$x\notin\ga_{A_iA_j}$. Therefore all the vertices of $P$, belonging to the same component  $C_l$ of the set $K\mmm \{x\}$, lie also in the same component of  $\hat\ga\mmm
\{x\}$, which implies $Ord(x,{\hat\ga})=Ord(x,K)$. \vse\\

 To evaluate the order $Ord(x,K)$ of the points $x\in K$, first we have to evaluate the order of the vertices   $A\in V_P$, which is related to the number of preimages  $n_A=\#\pi^{-1}(A)$ of the point $A$ in $\io$, and we evaluate it using measures $\te_A$ of solid angles at the vertices of $P$.
 
 Let $\te_A=\te(\Om(P,A))$ be the $(d-1)$-dimensional   measure of solid angle of $P$ at $A$,
   $\te_{max}=\max\{\te_A, A \in V_P\}$, and $\te_{min}=\min\{\te_A, A \in V_P\}$. 
   
For $t\in\rr$,  $\lceil t\rceil$ means  {\em $Ceil(t)$}, i.e. minimal integer, less or equal to $t$.

\begin{prop}\label{comp} Let $A \in V_P$.\\
 (a) If $\#\pi^{-1}(A )=1$, then there are $\bi\in\ia$, $A'\in V_P$, such that\\ $A=S_\bi(A')$ and $Ord(A,K)=Ord(A',\hat\ga)$; then
  $Ord(A ,K)\le n_P-1$;\\
 (b) If $n_A=\#\pi^{-1}(A )>1$,then there are $\bi_k\in\ia$, $A_k'\in V_P$, where $k=1,\ldots,n_A$, such that $A_k=S_{\bi_k}(A_k')$ and $Ord(A,K)=\sum\limits_{k=1}^{n_A} Ord(A_k',\hat\ga)$; then 
 $$Ord(A, K)\le
(n_P-1)\left(\left\lceil{\dfrac{\te_{A}}{\te_{min}}}\right\rceil-1\right) \le
(n_P-1)\left(\left\lceil{\dfrac{\te_{max}}{\te_{min}}}\right\rceil-1\right)\eqno{(1)}$$
 \end{prop}
\dok
Let $\#\pi^{-1}(A )=1$ and $\{C_l,l=1,...,k\}$   be some finite set of components of $K\mmm\{
A\}$. Since $\{A \}$ is the intersection of unique sequence  of polyhedra  $P_{j_1}\NI P_{j_1j_2}\NI...\NI
P_{j_1..j_s}..$, there is such $s$, that $\diam
P_{j_1..j_s}<\diam C_l$ for any $l=1,...,k$. Then, by Proposition \ref{squeeze}, each component $C_l$ contains a vertex
of a polyhedron $P_{j_1..j_s}$, different from  $A $. Therefore $k\le
n_P-1$, and $Ord(A,K)\le n_P-1$. 

  Since $Ord(A,K)$ is finite, we have the right to suppose that  $k=Ord(A,K)$, and  $\{C_1,...,C_k\}$ is a complete set of components of  $K\setminus \{A\}$.\medskip

Let $\bj=j_1..j_s$ and $A =S_\bj(A')$. The sets $C_l\cap P_\bj, l=1,...,k$
are the components of $K_\bj\mmm\{A\}$. Since $(K\cap
P_\bj)\mmm\{A \}=S_\bj(K\mmm\{A'\})$, the set $K\mmm \{A'\}$ consists of  $k$ 
components  $C_l'$, such that $S_\bj(C_l')=C_l\cap
P_\bj$. Since each component  $C_l'$ contains vertices of $P$, by Proposition \ref{xintree}(d), $Ord(A',{\hat\ga})=Ord(A',K)=Ord(A,K)\le
n_P-1$. \medskip

Suppose that $n_A=\#\pi^{-1}(A)>1$. By Theorem \ref{refsys}(g) there is a family $\{\Om_1,\ldots,\Om_{n_A}\}$  of disjoint solid angles with the same vertex $ A  $, and of respective polyhedra $P_{\bj_k}\ni A$, such that $P_{\bj_k}\IN\Om_k$ and $\Om(P_{\bj_k},A)=\Om_k$. 

Let $A_k\in V_P$  and $ S_{\bj_k}(A_k)=A$. Keeping in mind that $\#\pi^{-1}(A_k)=1$ and following the argument of the part (a)we can choose such $\bj_k$ and $A_k'$ that $Ord(A',K)=Ord(A'_k,\hat\ga)$; 
therefore $Ord(A,K_{\bj_k})=Ord(A_k,K)\le n_P-1$ and $Ord(A,K)\le n_A (n_P-1)$.  Taking into account the inequality 
$n_A\le \left\lceil{\dfrac{\te_A}{\te_{min}}}\right\rceil-1\le
\left\lceil{\dfrac{\te_{max}}{\te_{min}}}\right\rceil-1$, we get the inequality (1).\vse

\begin{thm}\label{order}

(i)\ \ $CP(K)\IN G_\eS(\hat\ga)$;\\
 (ii)\ \ If $y\notin G_\eS(V_P)$, then there are $\bj\in\ia$, $x\in CP(\hat\ga)$,\\ such that $y=S_\bj(x)$ and  $Ord(y,K)=Ord(x,{\hat\ga})\le n_P$.\\ 
 (iii)\ \ If $y\in G_\eS(V_P)$, then there are multiindices $\bj_k, k=1,..,s$ and vertices $A'_1,...,A'_s$, such that for any $k$, $S_{\bj_k}(A'_k)=y$, and any $l\neq k$,
$S_{\bj_k}(P)\cap S_{\bj_l}(P)=\{y\}$;\\ in this case, $
Ord(y,K)=\sum\limits_{k=1}^s
Ord(A'_k,{\hat\ga})\le(n_P-1)\left(\left\lceil{\dfrac{\te_{F}}{\te_{min}}}\right\rceil-1\right)$, where $\te_{F}$ is the measure of full angle in $\rd$.
\end{thm}

{\bf Proof.} (ii) Let $\{C_1,...,C_k\}$ be some set of
components of $K\setminus \{y\}$, and  $\rho=\min\limits_{l=1,...,k}\diam(C_l)$. Suppose $\bj\in I^*$  is such that
 $y\in P_\bj$ and $\diam(P_\bj)<\rho$.
 
 By
Proposition \ref{squeeze}, for any  $l$,   $C_l\cap S_\bj(V_P)\neq\0$, therefore $k\le n_P$. Thus, $Ord(y,K)\le n_P$.
  
 So we can suppose
that  $k=Ord(y,K)$
 and $\{C_1,...,C_k\}$ is the set of all components of
$K\setminus \{y\}$. 

Let $x=S_\bj^{-1}(y)$. Then the sets $C_l'=S_\bj^{-1}(C_l\cap P_\bj), l=1,\ldots,k,$ form a full set of components of $K\mmm \{x\}$, while for any  $l$, $C_l'\cap V_P\neq\0$. Then, by Proposition
\ref{xintree},
$Ord(x,{\hat\ga})=Ord(x,K)=Ord(y,K)\le n_P$.

(iii) Let $n_y=\#\pi^{-1}(y)$. By  Theorem \ref{refsys}(g) there is a family $\{\Om_1,\ldots,\Om_{n_y}\}$  of disjoint solid  angles with vertex $ y  $, and of respective  polyhedra $P_{\bj_k}\ni y$, such that $P_{\bj_k}\IN\Om_k$ and $\Om(P_{\bj_k},y)=\Om_k$. 

Using the argument of Proposition \ref{comp}(b), we obtain that $Ord(y,K)\le n_y(n_P-1)$ and therefore, choosing the polyhedra $P_{\bj_k}$ of sufficiently small diameter, we obtain that for any $k$, $y\in S_{\bj_k}(\hat\ga)$, $Ord(y,K_{\bj_k})= Ord(y, S_{\bj_k}(\hat\ga)$.
This gives the estimate $$Ord(y,K)\le(n_P-1)\left(\left\lceil{\dfrac{\te_{F}}{\te_{min}}}\right\rceil-1\right)$$
(i) follows from (ii) and (iii).\vse

\begin{thm}\label{equal}
Let $(P,\eS)$ be a contractible $P$-polyhedral   system and $K$ be its attractor.
(i) $\dim_H(CP(K))=\dim_H(\hat\ga)\le \dim_H EP(K)=\dim_H(K)$;
(ii) $\dim_H(CP(K))=\dim_H(K)$ iff $K$ is a Jordan arc.
\end{thm}
\dok Since $CP(K)=G_\eS(\hat\ga)$,  $\dim_H(CP(K))=\dim_H(\hat\ga)$.
If $K$ is not a Jordan arc, the set   $EP(K)$ is infinite \cite[Theorem 5.2]{Hata} and contains a point  $x\notin\hat\ga$. Let $\ep<d(x,\hat\ga)/2$. Take such $n$ that for any $\bj\in I^n$, 
$\diam(P_\bj)<\ep$. Then the set $\eJ=\{\bj\in I^n: P_\bj\cap\hat\ga\neq\0\}$ is a proper subset of $I^n$, because $x\notin P_\bj$ for any $\bj\in\eJ$. Let $\eS'=\{P_\bj,\bj\in\eJ\}$ and $K'$ be the attractor of the system  $\eS'$. Since the sets $\{S_\bj,\bj\in\eJ\}$ cover $\hat\ga$, 
$K'\NI\hat\ga$. At the same time, the similarity dimension  $\dim_s(\eS')$ of the system $\eS'$ is strictly smaller than that of $\eS^{(n)}$ which is equal to $\dim_s(\eS)=\dim_H(K)$ in its turn. Therefore, 
$\dim_H(\hat\ga)\le\dim_H(K')<\dim_H(K)$. Since $\\EP(K)=K\mmm CP(K)$, $\dim_H(EP(K))=\dim_H(K)$.

\vse

{\vspace{.1in}
\noindent Mary Samuel\\
Department of Mathematics\\
Bharata Mata College, Kochi, India\\
marysamuel2000@gmail.com

\vspace{.1in}\noindent Andrei Tetenov\\
Gorno-Altaisk State University 
and\\
Novosibirsk State University, Russia\\
atet@mail.ru\\

\vspace{.1in}\noindent Dmitry Vaulin\\
Gorno-Altaisk State University\\
Gorno-Altaisk, Russia\\
d\_warrant@mail.ru }

\end{document}